\documentclass[10pt,a4paper]{article}

\usepackage{amsmath}
\usepackage{amsfonts}

\newtheorem{theorem}{Theorem}[section]
\newtheorem{lemma}{Lemma}[section]
\newtheorem{corollary}{Corollary}[section]

\newenvironment{proof}[1][Proof]{\noindent \textit{#1.} }{\hfill $\Box$}

\renewcommand{\P}{\mathbb{P}}
\newcommand{\E}{\mathbb{E}}
\newcommand{\F}{\mathcal{F}}
\newcommand{\EE}{\mathcal{E}}
\newcommand{\A}{A}

\title{Weighted moments of the limit of a Branching Process in a Random
Environment}
\author{
Xingang Liang$^{\textrm{a,b}}$ \;\; and \;\; Quansheng Liu$^{\textrm{a,b,}}$
\thanks{
Corresponding author at: LMAM, Universit\'e de Bretagne-Sud, Campus de Tohannic, BP 573, 56017 Vannes, France. Tel.: +33 2 9701 7140; fax: +33 2 9701 7175.
\newline  {\em E-mail addresses:} xingang.liang@univ-ubs.fr (X. Liang), quansheng.liu@univ-ubs.fr (Q. Liu).
}
\\
{\small
$^{\textrm{a}}$LMAM, Universit\'e de Bretagne-Sud,
Campus de Tohannic,
 BP 573, 56017 Vannes, France}\\
{\small $^{\textrm{b}}$Universit\'e Europ\'eenne de Bretagne, France}
}
\date{}

\begin{document}

\maketitle

\begin{abstract}
\noindent Let $(Z_n)$ be a supercritical branching process in a random environment $%
\zeta$, and $W$ be the limit of the normalized population size $Z_n/\mathbb{E%
}(Z_n|\zeta)$. We show necessary and sufficient conditions for the existence
of weighted moments of $W$ of the form $\E W^{\alpha}\ell(W)$, where
$\alpha\geq 1$, $\ell$ is a positive function slowly varying at $\infty$. In the Galton-Watson case, the results improve those of Bingham and Doney (1974).


\end{abstract}

\section{Introduction and main results}

\setcounter{equation}{0}

For a Galton-Watson process  $(Z_n)$  with  offspring mean $m  = \E Z_1  \in (1, \infty)$,
the moments of $W = \lim Z_n/m^n$  have been studied by many authors: see for example
Harris (1963), Athreya and Ney (1972),  Bingham and Doney (1974), Alsmeyer and R\"osler (2004).  Bingham and Doney (1974) established very interesting comparison theorems between $W$ and $Z_1$ by considering weighted moments of the form $\E W^\alpha \ell(W)$, where $\alpha \geq 1$ and $\ell$ is a positive function slowly varying at $\infty$.   In particular, for $\alpha >1$, they showed that $\E W^\alpha \ell(W) < \infty$ if and only if $\E Z_1^\alpha \ell(Z_1) < \infty$, whenever
$\alpha$ is not an integer or $ \ell(x) = \int_1^x \ell_0(x)/x dx$ for some function $\ell_0$ slowly varying at $\infty$.  Alsmeyer and R\"osler (2004) showed that  the additional condition on $\ell$ can be removed if $\alpha$ is not a power of $2$.
In this paper, we shall show that this condition can always be removed. However,  our main objective is  to prove similar results for a branching process in a random environment.

Let $\zeta =(\zeta _{0},\zeta _{1},\ldots )$ be a sequence of independent
and identically distributed (i.i.d.) random variables, taking values in some space $\Theta$, whose realization corresponds to a sequence of
probability distributions on $\mathbb{N}$:
\begin{equation}
p{(\zeta _{n})}=\{p_{i}(\zeta _{n}):i\geq0\}, \text{ where } p_{i}(\zeta
_{n})\geq 0,\ \sum_{i=0}^{\infty }p_{i}(\zeta _{n})=1.
\end{equation}
A branching process $(Z_{n})_{n\geq 0}$ in the random environment $\zeta $
(BPRE) is a family of time-inhomogeneous branching processes (see e.g. \cite{Athreya71a,Athreya71b,Athreya72}): given the
environment $\zeta$, the process $(Z_{n})_{n\geq 0}$ acts as a Galton-Watson
process in varying environments with offspring distributions $p({\zeta _{n}}%
)$ for particles in $n$th generation, $n\ge0$. By definition,
\begin{equation}
Z_{0}=1\qquad {\rm and}\qquad Z_{n+1}=\sum_{u\in T_n}X_{u}\quad \mathrm{%
for}\quad n\geq 0,
\end{equation}
where conditioned on $\zeta $, $\{X_{u}:|u|=n\}$ are integer-valued random
variables with common distribution $p({\zeta _{n}});$ all the random
variables $X_u$, indexed by finite sequences of integers $u$, are
conditionally independent of each other. Here $T_n$ denotes the set of all
individuals of generation $n$, denoted by sequences $u$ of positive integers of
length $|u|=n$: as usual, the initial particle is denoted by the empty sequence \o\ (of length $0$); if $u\in T_n$, then $ui\in T_{n+1}$ if and only if $1\le i\le X_u$. The classical Galton-Watson process corresponds to the case where all $\zeta_n$ are the same constant.

Let $(\Gamma ,\P _{\zeta })$ be the probability space under which the
process is defined when the environment $\zeta $ is given. Therefore under $\P_{\zeta}$, the random variables $X_u$ are independent of each other, and have the common law $p(\zeta_n)$ if $|u|=n$. The probability $%
\P _{\zeta }$ is usually called \emph{quenched law}. The total probability
space can be formulated as the product space $(\Theta ^{{\mathbb{N}}}\times
\Gamma ,\P )$, where $\P =\P _{\zeta }\otimes \tau $ in the sense that for
all measurable and positive function $g$, we have
\begin{equation*}
\int gdP=\int \int g(\zeta ,y)d\P _{\zeta }(y)d\tau (\zeta ),
\end{equation*}%
where $\tau $ is the law of the environment $\zeta $. The total probability $\P $ is called \emph{%
annealed law}. The quenched law $\P _{\zeta }$ may be considered to be the
conditional probability of the annealed law $\P $ given $\zeta $. The
expectation with respect to $\P _{\zeta }$ (resp. $\P $ ) will be denoted $%
\E_{\zeta }$ (resp. $\E$).

For $n\ge0$, write
\begin{equation}
m_{n}=\sum_{i=0}^{\infty }ip_{i}(\zeta _{n}),\quad \Pi _{0}=1 \quad {\rm and}\quad \Pi _{n}=m_{0}\cdots m_{n-1}\ {\rm if }\ n\geq 1.
\end{equation}%
Then $E_{\zeta}X_u=m_n$ if $|u|=n$, and $E_{\zeta }Z_{n}=\Pi _{n}$ for each $%
n $.

We consider the supercritical case where
\begin{equation*}
\E\ln m_{0}\in(0,\infty].
\end{equation*}
It is well-known that under $\P _{\zeta }$,
\begin{equation*}
W_{n}=\frac{Z_{n}}{\Pi _{n}}\qquad (n\geq 0)
\end{equation*}%
forms a nonnegative martingale with respect to the filtration
\begin{equation*}
\mathcal{E}_{0}=\{\emptyset ,\Omega \}\qquad{\rm and}\qquad \mathcal{E}
_{n}=\sigma \{\zeta ,X_{u}:\ |u|<n\}\quad {\rm for}\ n\ge1.
\end{equation*}
It follows that $(W_n,\EE_n)$ is also a martingale under $\P$.
Let
\begin{equation}
W=\lim_{n\rightarrow\infty} W_{n}\qquad {\rm and}\qquad W^{\ast }:=\sup_{n\geq 0}W_{n},
\end{equation}%
where the limit exists a.s. by the martingale convergence theorem, and $\E W\leq 1$ by Fatou's lemma.

 We are interested in asymptotic properties of $W.$ Recall that in %
\cite{Liu01b}, Guivarc'h and Liu gave a necessary and sufficient condition
for the existence of moments of $W$ of order $\alpha >1$:
\vskip0.2cm

\begin{lemma}
\label{aprop1} \textrm{(\cite[Theorem 3]{Liu01b})} Let $(Z_{n})$ be a
supercritical branching process in an i.i.d. random environment. Let $\alpha
>1$. Then $0<\E W^{\alpha }<\infty $ if and only if  $\E m_{0}^{-(\alpha -1)}<1$ and $\E W_{1}^{\alpha }<\infty $.
\end{lemma}

This result suggests that under a moment condition on $m_{0},$ $W_{1}
$ and $W$ have the same asymptotic properties. In the following, we shall
establish comparison theorems between  weighted moments of $W_{1}$ and $W$.

Recall that  a positive and measurable function  $\ell $ defined on $[0,\infty)$ is called
slowly varying at $\infty $ if $\lim\limits_{x\rightarrow \infty
}\ell(\lambda x)/\ell(x)=1$ for all $\lambda>0$. (Throughout this paper, the term "positive" is used in the wide sense.) By the representation theorem (see \cite[Theorem1.3.1]{Bingham87}),  any slowly varying function $\ell $ is of the form
\begin{equation}\label{ell}
\ell (x)=c(x)\exp \left(\int_{a_0}^{x}\epsilon (t)dt/t\right),\qquad x> a_0,
\end{equation}%
where $a_0\ge0$, $c(\cdot)$ is measurable with $c(x)\rightarrow c$ for some constant $c\in(0,\infty)$, and $\epsilon
(x)\rightarrow 0$, as $x\rightarrow\infty$. The value of $a_0$ and those of $\ell(x)$ on $[0,a_0]$ will not be important; we always assume that $\ell$ is bounded on compact sets. For convenience, we often take $a_0=1$.

 We search for conditions under which $W$ has weighted moments of
the form $\E W^{\alpha }\ell (W),$ where $\alpha \geq 1,$  $\ell \geq
0$ is a function slowly varying at $\infty .\ $Notice that the function $%
c(x)$ in the representation of $\ell (x)$ has no influence on the finiteness
of the moments, so that  we can suppose without loss of generality that $%
c(x)=1.$

We first consider the case where $\alpha >1.$ As usual, for a set $A$, we
write Int$A$ for its interior.


\begin{theorem}
\label{athm1} Let $\alpha \in \mathrm{Int}\{a>1:\E m_{0}^{1-a}<1\}$
and $\ell :$ $[0,\infty )\mapsto \lbrack 0,\infty )$ be a function slowly
varying at $\infty $. 
Then the following assertions are equivalent:

\begin{itemize}
\item[\textrm{(a)}] $\E W_1^{\alpha}\ell(W_{1})<\infty $%
\negthinspace\ ;
\item[\textrm{(b)}] $\E W^{*\alpha}\ell(W^{*})<\infty $%
\negthinspace\ ;
\item[\textrm{(c)}] $0<\E W^{\alpha }\ell (W)<\infty $\negthinspace\
. \
\end{itemize}
\end{theorem}


 The result is sharp even for the classical Galton-Watson process (where $\zeta _{n}$ are the
same constant): in this case, it improves the corresponding result of Bingham and
Doney (1974) in the sense that they needed an additional assumption on $\ell
$ (which is equivalent to the hypothesis that $\epsilon (t)$ is positive and
slowly varying at $\infty $) when $\alpha $ is an integer. Alsmeyer and
R\"osler (2004) showed that this additional condition can be removed if $%
\alpha $ is not a dyadic power; our result shows that it can be removed for
all $\alpha $ and that the same conclusion holds in the random environment case.

We now consider the case where $\alpha =1,$ where the situation is different
as already shown by Bingham and Doney (1974) in the Galton-Watson case.

For a measurable function $\ell :[0,\infty )\mapsto [0,\infty )$, we set
\begin{equation}\label{hatell}
\hat{\ell}(x)=%
\begin{cases}
\int_{1}^{x}\frac{\ell (t)}{t}dt & \ \mathrm{if}\ x>1; \\
0 & \ \mathrm{if}\ x\leq 1.%
\end{cases}%
\end{equation}

We essentially deal with the case where $\ell $ is concave, which covers the case of slowly varying functions considered by Bingham and Doney (1974). (cf. Corollary 1.1 below)

\begin{theorem}
\label{athm2} Let $\ell $ be a positive and concave function defined on $%
[a_{0},\infty )$ for some $a_{0}\geq 0$. If $\E m_{0}^{-1}<1$ and $%
\E W_{1}\hat{\ell}(W_{1})<\infty $, then
\begin{equation*}
\E W^{\ast }\ell (W^{\ast })<\infty \qquad {\rm and}\qquad \E W\ell
(W)<\infty \,.
\end{equation*}
Moreover, in the case where $\ell $ is also slowly varying at $\infty ,$ the
moment condition $\E m_{0}^{-1}<1$ can be relaxed to $\E%
m_{0}^{-\delta _{0}}<\infty $ for some $\delta _{0}>0.$
\end{theorem}

As a corollary, we obtain:

\begin{corollary}
\label{acor1} Let $\ell :$ $[0,\infty )\mapsto $ $[0,\infty )$ be
nondecreasing and slowly varying at $\infty $, such that $\ell
(x)=\int_{1}^{x}\ell _{0}(t)dt/t$ for some function $\ell _{0}\ge0$
slowly varying at $\infty $. Assume that $\E m_{0}^{-\delta
_{0}}<\infty $ for some $\delta _{0}>0.$ If $\E W_{1}\hat{\ell}%
(W_{1})<\infty $, then
\begin{equation*}
\E W^{\ast }\ell (W^{\ast })<\infty \qquad {\rm and}\qquad \E W\ell
(W)<\infty \,.
\end{equation*}
\end{corollary}

 Corollary \ref{acor1} extends the sufficiency of Theorem 7 of
Bingham and Doney (1974)  where the classical Galton-Watson
process was considered. (See also Corollary 2.3 of Alsmeyer and R\"osler
(2004).)

Notice that $\ell (x)=\int_{1}^{x}\ell _{0}(t)dt/t$ for some
function $\ell _{0}$ slowly varying at $\infty $ if and only if the
function $x\ell ^{\prime }(x)$ is slowly varying at $\infty $; when $\ell $
is of the canonical form $\ell (x)=\exp (\int_{1}^{x}\epsilon (t)d%
t/t)$ with $\epsilon (t)\rightarrow 0,$ this is exactly the case where $%
\epsilon (t)$ is slowly varying at $\infty .$

Corollary \ref{acor1} is a direct consequence of Theorem \ref{athm2}. To see this, we can suppose that $\ell $ is of the form (\ref{ell})
with $c(x)=1.$ Hence $x\ell ^{\prime }(x)=\ell (x)\epsilon (x)$. Let $\psi
(x)=\inf \{\ell ^{\prime }(t):1\leq t\leq x\}$. Since $x\ell ^{\prime }(x)$ is
slowly varying, we have $\ell ^{\prime }(x)\sim \psi(x) $ (see \cite[Theorem 1.5.3]{Bingham87}), where $\psi $ is positive and nonincreasing; so $\ell
(x)\asymp \ell _{1}(x):=\int_{0}^{x}\psi (t)dt$, and $\ell _{1}$ is
a positive and concave function. Here, as usual, we write
\begin{equation}
f(x)\asymp g(x) \qquad{\rm if}\quad
0<\liminf\limits_{x\rightarrow \infty }\frac{f(x)}{g(x)}\leq
\limsup\limits_{x\rightarrow \infty }\frac{f(x)}{g(x)}<\infty\,,
\end{equation}
and $f(x)\sim g(x)$
if $\lim\limits_{x\rightarrow \infty }f(x)/g(x)=1$. Therefore we can apply Theorem \ref{athm2}
to $\ell _{1}$ to obtain the conclusion of the corollary.

By the same method, we can consider some slightly different classes of
functions. For example, we can show the following result similar to Theorem
1.1 of Alsmeyer and R\"osler (2004) where the Galton-Watson case was considered.

\begin{corollary}
\label{acor2} Let $\phi$ be positive and convex on $[0,\infty)$ with
positive concave derivative $\phi^{\prime }$ on $(0,\infty)$. Define
\begin{gather*}
\tilde{\phi}(x)\!=\left\{
\begin{array}{ll}
\int_1^x\frac{\phi^{\prime }(t)}{t}dt & \ \mathrm{if} \ x>1; \\
0 & \ \mathrm{if} \ x\le 1.%
\end{array}%
\right.
\end{gather*}
If $\E m_0^{-1}<1$ and $\E W_1\tilde{\phi}(W_1)<\infty $,
then
\begin{gather*}
\E\phi(W^*)<\infty\qquad {\rm and}\qquad \E\phi(W)<\infty\,.
\end{gather*}
\end{corollary}


The argument in the proof of Theorem \textsc{\ref{athm2} }can also be used
to study the integrability of $W^{\ast }$ and thus the non-degeneration of $W$. As usual, we write $\ln^-x=\max{(0,-\ln x)}$.

\begin{theorem}\label{athm3}
Assume that $\mathbb{E(}\ln^{-}m_{0})^{2}<\infty$. If $\E W_{1}\ln ^{+}W_{1}<\infty$, then $\E W^{\ast }<\infty $.
\end{theorem}

 Notice that $\E W^{\ast }<\infty $ implies $EW=1$ by the
dominated convergence theorem.\ Therefore Theorem \ref{athm3} implies the
classical theorem (the sufficiency) of Kesten-Stigum (1966) on the Galton-Watson process. It
gives a new proof of the corresponding result of Athreya and Karlin (1971b)
(see also Tanny (1988)) for a branching process in a random environment, under
the extra condition that $\E(\ln ^{-}m_{0})^{2}<\infty$. (Notice that the supercritical condition $\E \ln m_0>0$ implies $\E \ln^-m_0<\infty$.) Although
we need this extra condition, the conclusion that $\E W^{\ast }<\infty $
may be useful in applications; we do not know whether this conclusion is equivalent
to $EW=1.$ (It is known (see \cite{Tanny88})  that the condition $\E%
W_{1}\ln ^{+}W_{1}<\infty $ is equivalent to $EW=1$; in the Galton-Watson case, it is also known that this condition is also equivalent to $\E W^*<\infty$. But we do not know whether the same conclusion remains true for the random environment case.)

In the Galton-Watson case, Alsmeyer and R\"osler (2004) used a similar argument (also based on convex inequalities for martingales) to show the non-degeneration of $W$. But our approch is more direct, as we do not use their Lemma 4.5. 

\vskip 2mm
The rest of the paper is organized as follows. In Section 2, we establish key inequalities based on convex inequalities on martingales. In Section 3, we give corrected versions of regularly varying functions in order to use the key inequalities of Section 2. Theorem \ref{athm1} is proved in Section 4, while Theorems \ref{athm2} and \ref{athm3} are proved in Sections 5 and 6, respectively.

\vskip 2mm
In enclosing the introduction, we mention that the argument of this paper can be adapted to weighted branching processes, thus enables us to improve the results of Bingham and Doney (1975) for Crump-Mode and Jirina processes, those of Alsmeyer and Kuhlbusch (2009) for branching random walks, and to extend their results to the random environment case (including the weighted branching processes considered by Kuhlbusch(2004)). This will be done in the forth coming paper \cite{Liang10}.

\section{Key Inequalities}

\setcounter{equation}{0}

In this section, we show key inequalities that will be used for the proof of main theorems. As in Alsmeyer and R\"osler (2004), our argument is based on convex inequalities on martingales.

We first introduce some notations. For a finite sequence $u\in\bigcup_{n=0}^{\infty}{\mathbb{N}^*}^n$ (${\mathbb{N}^*}^0=\{\emptyset\}$ by convention), set $\tilde{X}_u\!=\frac{X_u}{m_{|u|}}-1$. For $n\ge1$, write
\begin{equation}
D_n\!=W_n-W_{n-1}=\frac{1}{\Pi_{n-1}}\sum_{|u|=n-1}\tilde{X}_u.
\end{equation}
Then $W^*=\sup_{n\ge0}W_n$ can be written as
\begin{equation*}
W^*=1+\sup_{n\ge1}(D_1+\ldots+D_n).
\end{equation*}
For convenience, let $\tilde{X}_n=\tilde{X}_{u_0|n}$ where $u_0\in\mathbb{N}^*\times\mathbb{N}^*\times\cdots$ is a fixed infinite sequence, $u_0|n$ denotes the restriction to the first $n$ terms of $u_0$.

Define
\begin{equation}
\mathcal{F}_0=\{\emptyset,\Omega\}\qquad {\rm and} \qquad \mathcal{F}%
_n=\sigma\{\zeta_k,\tilde{X}_{u}: k<n,|u|<n\}\quad {\rm for}\ n\ge1.
\end{equation}
Then $(W_n,\mathcal{F}_n)_{n\ge0}$ also forms a nonnegative
martingale under $\P$, as
\begin{equation*}
\E(W_n|\mathcal{F}_{n-1})=\E(\E(W_n|\mathcal{E}%
_{n-1})|\mathcal{F}_{n-1})=\E(W_{n-1}|\mathcal{F}_{n-1})=W_{n-1}.
\end{equation*}
For technical reasons, we will use the martingale $(W_n,\F_n)$, rather than the more frequently used one $(W_n, \mathcal{E}_n)$. We will explain this after the proof of Theorem \ref{athm1}. For convenience,  we shall write for $n\ge0$,
\begin{equation}
\P_n(\cdot)=\P(\cdot|\F_n)\qquad {\rm and} \qquad \E_n(\cdot)=\E(\cdot|\F_n).
\end{equation}
The letter $C$ will always denote a finite and positive constant which may differ from line to line. The terms ''increasing'' and ''decreasing'' will be used in the wide sense.


\begin{theorem}\label{alm1}
Let $\phi$ be convex and increasing with $\phi(0)=0$ and $%
\phi(2x)\le c\phi(x)$ for some constant $c\in(0,\infty)$ and all $x>0$. Let $\beta\in(1,2]$.

\begin{itemize}
\item[\textrm{{(i)}}] If the function $x\mapsto\phi(x^{1/\beta})$ is convex
and $\E|\tilde{X} _1|^{\beta}<\infty$, then writting $A=\sum_{n=1}^\infty\frac{1}{ \Pi_{n-1}^{\beta-1}}$, we have
\begin{eqnarray}  \label{aeqlm11}
\E\phi(W^*-1)&\leq& C\sum_{n=1}^{\infty}
\left(\E\left(\frac{%
1}{A\Pi_{n-1}^{\beta-1}}\phi(A^{1/\beta} W_{n-1}^{1/\beta})\right)\right.
\notag \\
&&\qquad\qquad \left.+\E \phi\left(\frac{|\tilde{X}_{n-1}|}{\Pi_{n-1}^{(\beta-1)/%
\beta}}\cdot W_{n-1}^{1/\beta}\right)\right),
\end{eqnarray}
where $C=C(\phi,\beta)>0$ is a constant depending only on $\phi$ and $\beta$.
\item[\textrm{{(ii)}}] If the function $x\mapsto\phi(x^{1/\beta})$ is
concave, then
\begin{equation}  \label{aeqlm12}
\E\phi(W^*-1)\leq C\sum_{n=1}^\infty\E\Pi_{n-1}\phi\left(%
\frac{|\tilde{X}_{n-1}|}{\Pi_{n-1}}\right),
\end{equation}
where $C=C(\phi,\beta)>0$ is a constant depending only on $\phi$ and $\beta$.
\end{itemize}
\end{theorem}

\begin{proof}
(i) By the Burkholder-Davis-Gundy (BDG) inequality (see \cite{Chow95}),
\begin{equation}\label{aeqlm13}
\E\phi(W^*-1)\le B\left(\E\phi\left(\left(\sum_{n=1}^{\infty}\E_{n-1}|D_n|^{\beta}\right)^{\frac{1}{\beta}}\right)+\sum_{n=1}^{\infty}\E\phi(|D_n|)\right),
\end{equation}
where $B>0$ is a constant depending only on $\phi$ and $\beta$.

Let $\tilde{X}(1),\ldots,\tilde{X}(Z_{n-1})$ be an enumeration of $\{\tilde{X}_u:u\in T_{n-1}\}$. By the fact that $\E_{\zeta}\tilde{X}(k)=0$ and the independence of $\{\tilde{X}_u\}$ under $\P_{\zeta}$, it can be easily seen that, under $P_{n-1}$, $\{\tilde{X}(1),\ldots,\tilde{X}(Z_{n-1})\}$ is a sequence of martingale differences with respect to the nature filtration
\begin{equation}
\tilde{\F}_k=\sigma\{\zeta_l, X_u:l<n-1, |u|<n-1, \tilde{X}(1),\ldots,\tilde{X}(k)\},\qquad k\ge 1.
\end{equation}
To this martingale difference sequence, using the BDG-inequality, we obtain
\begin{eqnarray}\label{aeqlm14}
\E_{n-1}|D_n|^{\beta}&=&\E_{n-1}\left|\frac{\sum_{|u|=n-1}\tilde{X}_u}{\Pi_{n-1}}\right|^{\beta}\nonumber\\
&\le& B\E_{n-1}\sum_{|u|=n-1}\frac{|\tilde{X}_u|^{\beta}}{\Pi_{n-1}^{\beta}}\nonumber\\
&=&B\frac{Z_{n-1}}{\Pi_{n-1}^{\beta}}\cdot\E_{n-1}|\tilde{X}_{n-1}|^{\beta}\nonumber\\
&=&C\frac{W_{n-1}}{\Pi_{n-1}^{\beta-1}},
\end{eqnarray}
where $C=B \E|\tilde{X}_1|^{\beta}<\infty$.
Since $\phi(x^{1/\beta})$ is convex and $\sum_{n=1}^{\infty}\frac{1}{\A\Pi_{n-1}^{\beta-1}}=1$, it follows that
\begin{eqnarray}\label{aeqlm15}
\E\phi\left(\left(\sum_{n=1}^{\infty}\E_{n-1}|D_n|^{\beta}\right)^{\frac{1}{\beta}}\right)&
\le&  \E\phi\left(\left(\sum_{n=1}^{\infty}\frac{C}{\A\Pi_{n-1}^{\beta-1}}\cdot \A W_{n-1}\right)^{\frac{1}{\beta}}\right)\nonumber\\
&\le& C\E\sum_{n=1}^{\infty}\frac{1}{\A\Pi_{n-1}^{\beta-1}}\cdot\phi\left(\A^{1/\beta}W_{n-1}^{1/\beta}\right).
\end{eqnarray}

For the second part of ($\ref{aeqlm13}$), again by the BDG-inequality and the convexity of $\phi(x^{1/\beta})$, we have
\begin{eqnarray}\label{aeqlm16}
\E_{n-1}\phi(|D_n|)&\le& B\E_{n-1}\phi\left(\left(\sum_{|u|=n-1}\frac{|\tilde{X}_u|^{\beta}}{\Pi_{n-1}^{\beta}}\right)^{\frac{1}{\beta}}\right)\nonumber\\
&\le& B\E_{n-1}\sum_{|u|=n-1}\frac{1}{Z_{n-1}}\phi\left(\frac{|\tilde{X}_u|}{\Pi_{n-1}}\cdot Z_{n-1}^{1/\beta}\right)\nonumber\\
&=& B\E_{n-1}\phi\left(\frac{|\tilde{X}_{n-1}|}{\Pi_{n-1}^{(\beta-1)/\beta}}\cdot W_{n-1}^{1/\beta}\right).
\end{eqnarray}
Therefore
\begin{equation}
\E\phi(|D_n|)\le B\E \phi\left(\frac{|\tilde{X}_{n-1}|}{\Pi_{n-1}^{(\beta-1)/\beta}}W_{n-1}^{1/\beta}\right).
\end{equation}

(ii) By the BDG-inequality and the concavity of $\phi(x^{1/\beta})$ (which implies the subadditivity),
\begin{eqnarray}\label{aeqlm17}
\E\phi(W^*-1)
&\leq& B\E\phi\left(\left(\sum_{n\geq1}|D_n|^\beta\right)^{\frac{1}{\beta}}\right)\nonumber\\
&\le& B\sum_{n\geq1}\E\phi(|D_n|),
\end{eqnarray}
where $B>0$ is a constant depending only on $\phi$ and $\beta$.

Similarly to the proof in part (i), by the  BDG-inequality and the concavity of $\phi(x^{1/\beta})$,
\begin{eqnarray}\label{aeqlm18}
\E_{n-1}\phi(|D_n|)&=&\E_{n-1}\phi\left(\left|\frac{1}{\Pi_{n-1}}\sum_{|u|=n-1}\tilde{X}_u\right|\right)\nonumber\\
&\le& B \E_{n-1}\phi\left(\left(\sum_{|u|=n-1}\frac{|\tilde{X}_u|^{\beta}}{\Pi_{n-1}^{\beta}}\right)^{\frac{1}{\beta}}\right)\nonumber\\
&\le& B \E_{n-1}\sum_{|u|=n-1}\phi\left(\frac{|\tilde{X}_u|}{\Pi_{n-1}}\right).
\end{eqnarray}
By the identical distribution of $(\tilde{X}_u)_{|u|=n-1}$ and the independence between $(\tilde{X}_u)_{|u|=n-1}$ and $Z_{n-1}$ under $\P$, we have
\begin{equation}\label{aeqlm19}
\E\phi(|D_n|)\le B\E Z_{n-1}\phi\left(\frac{|\tilde{X}_{n-1}|}{\Pi_{n-1}}\right)=\E \Pi_{n-1}\phi\left(\frac{|\tilde{X}_{n-1}|}{\Pi_{n-1}}\right).
\end{equation}
Combining (\ref{aeqlm17}) and (\ref{aeqlm19}), we get (\ref{aeqlm12}).

\end{proof}

\section{Corrected versions of regularly varying functions}

\setcounter{equation}{0} In this section, we will give some corrected versions of regularly varying functions to have better
properties.


\begin{lemma}
\label{alm2} Let $\phi(x)=x^\alpha\ell(x)$, with $\alpha>1$, and $%
\ell(x)=\exp\left(\int_{a_0}^x\epsilon(u)du/u\right)$ $(x\geq
a_0\geq0)$ with $\epsilon(x)\rightarrow0$ $(x\rightarrow0)$. Then for each $%
\beta\in(1,2]$ with $\beta<\alpha$, there is a function $\phi_1(x)\geq0$
such that:
\begin{itemize}
\item[\textrm{(i)}] $\phi_1(x)\sim\phi(x)$;
\item[\textrm{(ii)}] $\phi_1(x)$ and $\phi_1(x^{1/\beta})$ are convex on $%
[0,\infty)$;
\item[\textrm{(iii)}] $\phi_1(x)=x^\alpha\ell_1(x)$, where $\ell_1(x)$ is
slowly varying at $\infty$ and $\ell_1(x)>0$ $\forall x\geq0$.
\end{itemize}
\end{lemma}

\begin{proof}
Fix $\beta\in(1,2]$ with $\beta<\alpha$. Notice that the derivative
\begin{equation*}
\phi'(x)=x^{\alpha-1}\ell(x)(\alpha+\epsilon(x))
\end{equation*}
behaves like $\alpha x^{\alpha-1}\ell(x)$ as $x\rightarrow\infty$. It is therefore natural to define
\begin{equation}\label{aeqlm21}
\phi_1(x)=\alpha\int_0^x t^{\alpha-1}\ell(t)dt, \qquad x>a,
\end{equation}
where $a\geq\max(1,a_0)$ is large enough such that $\forall\ x>a$, $\alpha-\beta+\epsilon(x)>0$, so that
\begin{equation}\label{aeqlm22}
\frac{d}{dx}(x^{\alpha-1}\ell(x))=x^{\alpha-2}\ell(x)(\alpha-1+\epsilon(x))>0 \qquad \forall x>a,
\end{equation}
and
\begin{equation}\label{aeqlm23}
\frac{d}{dx}\left(x^{\frac{\alpha}{\beta}-1}\ell(x^{\frac{1}{\beta}})\right)=x^{\frac{\alpha}{\beta}-2}\ell(x^{\frac{1}{\beta}})\left((\frac{\alpha}{\beta}-1)+\frac{\epsilon(x^{\frac{1}{\beta}})}{\beta}\right)>0 \quad \forall x>a^\beta.
\end{equation}
Therefore, $\phi_1(x)$ is convex on $(a,\infty)$ as $\phi_1'(x)=\alpha x^{\alpha-1}\ell(x)$ is increasing; and $\phi_1(x^{1/\beta})$ is convex on $[a^\beta,\infty)$ as
\begin{equation*}
\frac{d}{dx}\phi_1(x^{1/\beta})=\phi_1'(x^{1/\beta})\cdot\frac{1}{\beta}x^{\frac{1}{\beta}-1}=\frac{\alpha}{\beta}x^{\frac{\alpha}{\beta}-1}\ell(x^{1/\beta})\qquad (x>a^\beta)
\end{equation*}
is also increasing on $(a^\beta,\infty)$. Define
\begin{equation}\label{aeqlm24}
\phi_1(x)=x^\alpha\ell(a),\qquad\forall x\in[0,a].
\end{equation}
Then
\begin{equation}\label{aeqlm25}
\frac{d}{dx}\phi_1(x)=\alpha x^{\alpha-1}\ell(a)\qquad\forall x\in[0,a],
\end{equation}
and
\begin{eqnarray}\label{aeqlm26}
\frac{d}{dx}\phi_1(x^{1/\beta})&=&\frac{d}{dx}\left(x^{\alpha/\beta}\ell(a)\right)\nonumber\\
&=&\frac{\alpha}{\beta}x^{\frac{\alpha}{\beta}-1}\ell(a)\qquad\forall x\in[0,a^\beta].
\end{eqnarray}
It follows that both $\frac{d}{dx}\phi_1(x)$ and $\frac{d}{dx}\phi_1(x^{1/\beta})$ are increasing on $[0,\infty)$. Therefore both $\phi_1(x)$ and $\phi_1(x^{1/\beta})$ are convex on $[0,\infty)$. Moreover,
\begin{equation}\label{aeqlm27}
\lim_{x\rightarrow\infty}\frac{\phi_1(x)}{\phi(x)}=\lim_{x\rightarrow\infty}\frac{\phi_1'(x)}{\phi'(x)}=1,
\end{equation}
so that $\phi_1(x)=x^\alpha\ell_1(x)$ for some slowly varying function $\ell_1$. If $x> a$, then $\ell_1(x)>0$ as $\phi_1(x)>0$; if $x\leq a$, then $\ell_1(x)=\ell(a)>0$. Therefore, $\ell_1(x)>0$ $\forall x\geq0$.

\end{proof}

\begin{lemma}
\label{alm3} Let $\ell$ be a positive and increasing function on $[0,\infty)$, concave
on $(a_0,\infty)$ for some $a_0\geq0$. Then there is a convex increasing
function $\phi_1(x)\geq0$ such that:
\begin{itemize}
\item[\textrm{(i)}] $\phi_1(x)\asymp x\ell(x)$;
\item[\textrm{(ii)}] $\phi_1(2x)\le c\phi_1(x)$ for some constant $%
c\in(0,\infty)$ and all $x>0$;
\item[\textrm{(iii)}] $\phi_1(x^{1/2})$ is concave on $(0,\infty)$.
\end{itemize}
\end{lemma}

\begin{proof}
Let
\begin{equation}\label{aeqlm31}
\ell_1(x)=\left\{
\begin{array}{ll}
\ell'(a)x  \qquad&{\rm if}\quad x\in [0,a],\\
\ell(x)+c_0 &{\rm if}\quad x\in(a,\infty),
\end{array}\right.
\end{equation}
where $a>a_0>0$, $c_0=\ell'(a)a-\ell(a)$, and $\phi_1(x)=\int_0^x\ell_1(t)dt$. We will show that $\phi_1$ satisfies the associate properties.

First, $\phi_1$ is convex as $\phi_1'(x)=\ell_1(x)$ is increasing on $[0,\infty)$; $\phi_1$ is increasing as $\ell_1$ is positive on $[0,\infty)$.

Next, for $x>2a$, as $\ell$ is increasing, we have $\ell_1(t)\ge\ell_1(\frac{x}{2})=\ell(\frac{x}{2})+c_0$ if $t\in[\frac{x}{2},x]$, and $\ell_1(t)\le \ell(x)$ if $t\in[0,x]$; therefore
\begin{equation}\label{aeqlm32}
\frac{x}{2}\left(\ell(\frac{x}{2})+c_0\right)\leq\phi_1(x)\leq x\ell(x).
\end{equation}
By the concavity of $\ell$, for all $x>a$,
\begin{eqnarray}\label{aeqlm33}
\ell(2x)
&=&\ell(x)+2\int_{\frac{x}{2}}^x\ell'(2s)ds\nonumber\\
&\leq&\ell(x)+2\int_0^x\ell'(s)ds\nonumber\\
&\leq&3\ell(x).
\end{eqnarray}
(\ref{aeqlm32}) and (\ref{aeqlm33}) imply that $\phi_1(x)\asymp x\ell(x)$ and that there is  a constant $c\in(0,\infty)$ such that $\ell_1(2x)\leq c\ell_1(x)$ for all $x>0$.

Moreover, we can prove that $\phi_1(x^{1/2})$ is concave. In fact,
\begin{eqnarray}\label{aeqlm34}
\frac{d}{dx}\phi_1(x^{1/2})
&=&\phi_1'(x^{\frac{1}{2}})\cdot \frac{1}{2}x^{-\frac{1}{2}}\nonumber\\
&=&\frac{1}{2}\ell_1(x^{\frac{1}{2}})\cdot x^{-\frac{1}{2}}.
\end{eqnarray}
Notice that $\frac{\ell_1(t)}{t}$ is decreasing as $\ell_1$ is concave with $\ell_1(0)=0$; hence $\frac{d}{dx}\phi_1(x^{1/2})$ is decreasing, so that $\phi_1(x^{1/2})$ is concave.
\end{proof}

{}

\section{Proof of Theorem \protect\ref{athm1}}

\setcounter{equation}{0}

\begin{proof}[Proof of Theorem \ref{athm1}]
Let $\beta\in(1,2]$ with $\beta< \alpha$. Write $\phi(x)=x^{\alpha}\ell(x)$. By Lemma \ref{alm2}, we can assume that the functions $\phi(x)$ and  $\phi(x^{1/\beta})$ are convex on $[0,\infty)$, and $\ell(x)>0$ $\forall\ x\ge0$.

(i) We first show that (a) implies (b). By Lemma \ref{alm1}, we obtain
\begin{eqnarray}\label{aeqthm12}
\E\phi(W^*-1)
&\le& C\sum_{n=1}^{\infty}\left(\E\left(\frac{1}{\A\Pi_{n-1}^{\beta-1}}\phi\left(\A^{1/\beta}W_{n-1}^{1/\beta}\right)\right)\right.\nonumber\\
&&\left.\qquad\qquad+\E\phi\left(\frac{|\tilde{X}_{n-1}|}{\Pi_{n-1}^{(\beta-1)/\beta}}\cdot W_{n-1}^{1/\beta}\right)\right).
\end{eqnarray}
Notice that $\ell>0$ on any compact subset of $[0,\infty)$, so by Potter's Theorem (see \cite{Bingham87}), for $\delta>0$ which will be determined later, there exists $C=C(\ell,\delta)>1$ such that $\ell(x)\leq C\max(x^{\delta},x^{-\delta})$ for all $x>0$. Hence for the first part of (\ref{aeqthm12}), we have
\begin{eqnarray}\label{aeqthm13}
\E\left(\frac{1}{\A\Pi_{n-1}^{\beta-1}}\phi\left(\A^{1/\beta}W_{n-1}^{1/\beta}\right)\right)&=&\E\left(\Pi_{n-1}^{1-\beta}\A^{\frac{\alpha}{\beta}-1}W_{n-1}^{\frac{\alpha}{\beta}}\ell\left(\A^{1/\beta}W_{n-1}^{1/\beta}\right)\right)\nonumber\\
&\le& C(I_1^+(n)+I_1^-(n)),
\end{eqnarray}
where
\begin{eqnarray*}
I_1^+(n)=\E\Pi_{n-1}^{1-\beta}\A^{\frac{\alpha+\delta}{\beta}-1}W_{n-1}^{\frac{\alpha+\delta}{\beta}},\\
I_1^-(n)=\E\Pi_{n-1}^{1-\beta}\A^{\frac{\alpha-\delta}{\beta}-1}W_{n-1}^{\frac{\alpha-\delta}{\beta}}.
\end{eqnarray*}
Recall that $Z_{n-1}$ is an integer-valued random variable with $\E_\zeta Z_{n-1}=\Pi_{n-1}$. Choose $\delta>0$ small enough such that $\beta-1-2\delta>0$. Then by H\"older's inequality, we obtain
\begin{eqnarray}\label{aeqthm13+}
\E_\zeta Z_{n-1}^{\frac{\alpha+\delta}{\beta}}
&=&\E_\zeta\left(Z_{n-1}^{\frac{\alpha+\delta-(\beta-1)}{\beta}}\cdot Z_{n-1}^{\frac{\beta-1}{\beta}}\right)\nonumber\\
&\leq&(\E_\zeta Z_{n-1}^{\alpha+\delta-(\beta-1)})^{\frac{1}{\beta}}\cdot(\E_\zeta Z_{n-1})^{\frac{\beta-1}{\beta}}\nonumber\\
&\leq&\left(\E_\zeta Z_{n-1}^{\alpha+\delta-(\beta-1)+(\beta-1-2\delta)}\right)^{\frac{1}{\beta}}\cdot(\E_\zeta Z_{n-1})^{\frac{\beta-1}{\beta}}\nonumber\\
&=&\Pi_{n-1}^{\frac{\alpha+\beta-1-\delta}{\beta}}\left(\E_\zeta W_{n-1}^{\alpha-\delta}\right)^{\frac{1}{\beta}}.
\end{eqnarray}
Therefore,
\begin{eqnarray}\label{aeqthm14}
I_1^+(n)&=&\E\left(\Pi_{n-1}^{1-\beta}\A^{\frac{\alpha+\delta}{\beta}-1}\cdot\frac{Z_{n-1}^{\frac{\alpha+\delta}{\beta}}}{\Pi_{n-1}^{\frac{\alpha+\delta}{\beta}}}\right)\nonumber\\
&\le& \E\left(\Pi_{n-1}^{-\frac{(\beta-1)^2+2\delta}{\beta}}\A^{\frac{\alpha+\delta}{\beta}-1}\left(\E_{\zeta}W_{n-1}^{\alpha-\delta}\right)^{\frac{1}{\beta}}\right).
\end{eqnarray}
Using H\"older's inequality twice, we see that
\begin{eqnarray}\label{aeqthm15}
I_1^+(n)&\le& \left(\E W_{n-1}^{\alpha-\delta}\right)^{\frac{1}{\beta}}\cdot\left(\E\Pi_{n-1}^{-\frac{(\beta-1)^2+2\delta}{\beta-1}}\A^{\frac{\alpha+\delta-\beta}{\beta-1}}\right)^{\frac{\beta-1}{\beta}}\nonumber\\
&\le& \left(\E W_{n-1}^{\alpha-\delta}\right)^{\frac{1}{\beta}}\cdot\left(\E \Pi_{n-1}^{-\frac{(\beta-1)^2+2\delta}{\beta-1}p}\right)^{\frac{\beta-1}{p\beta}}\cdot\left(\E \A^{\frac{\alpha+\delta-\beta}{\beta-1}p^*}\right)^{\frac{\beta-1}{p^*\beta}},
\end{eqnarray}
where $p>1$, $p^*>1$ and $\frac{1}{p}+\frac{1}{p^*}=1$. By Potter's Theorem, there exists $C=C(\ell,\delta)>0$ such that $\ell(x)\geq Cx^{-\delta}$ for all $x\geq1$. This yields
\begin{equation}\label{aeqthm16}
\E|\tilde{X}|^{\alpha-\delta}\leq  C(1+\E\phi(|\tilde{X}|))<\infty\,.
\end{equation}
Since $\alpha\in{\rm Int}\{a>1:\E m_0^{1-a}<1\}$, there exists $\delta_0\in(0,1)$ such that
\begin{gather*}
\E m_0^{1-(\alpha+\delta_0)}<1.
\end{gather*}
Notice that the function $\rho(x)=\E m_0^{1-x}$ is convex with $\rho(1)=1$, so  $\rho(\alpha+\delta_0)<1$ implies $\rho(x)<1$ for all $1<x<\alpha+\delta_0$; in particular, $\rho(\alpha-\delta)<1$.
Hence, by  Lemma \ref{aprop1},
\begin{gather}\label{aeqthm18}
\sup_{n\geq1}\E W_{n-1}^{\alpha-\delta}<\infty\,.
\end{gather}
We choose $p=1+\frac{(\alpha+\delta-\beta)(\beta-1)}{(\beta-1)^2+2\delta}$ so that $p_1:=\frac{\alpha+\delta-\beta}{\beta-1}p^*=\frac{(\beta-1)^2+2\delta}{(\beta-1)^2}p$. As $p_1(\beta-1)\in(1,\alpha+\delta_0)$ when $\delta$ is small enough, we get $\E\Pi_{n-1}^{-p_1(\beta-1)}=a^{n-1}$ with $a=\E m_0^{-p_1(\beta-1)}<1$; moreover, by the triangular inequality for the norm $\|\cdot\|_{p_1}$ in $L^{p_1}$,
\begin{equation*}
\|\A\|_{p_1}\le \sum_{n=1}^{\infty}\|\Pi_{n-1}^{-(\beta-1)}\|_{p_1}=\sum_{n=1}^{\infty}a^{(n-1)/p_1}<\infty\,.
\end{equation*}
Therefore,
\begin{equation}
\sum_{n=1}^{\infty}I_1^+(n)<\infty\,.
\end{equation}
We use a similar argument to estimate $I_1^-(n)$. This time, instead of (\ref{aeqthm14}), we have
\begin{equation}
I_1^-(n)\le \E\left(\Pi_{n-1}^{1-\beta}A^{\frac{\alpha-\delta}{\beta}-1}(\E_{\zeta} W_{n-1}^{\alpha-\delta})^{\frac{1}{\beta}}\right).
\end{equation}
Proceeding in the same way as before, we obtain
\begin{equation}
\sum_{n=1}^{\infty}I_1^-(n)<\infty\,.
\end{equation}
Hence
\begin{equation}\label{aeqthm110}
\sum_{n=1}^{\infty}\E\left(\frac{1}{\A\Pi_{n-1}^{\beta-1}}\phi\left(\A^{1/\beta}W_{n-1}^{1/\beta}\right)\right)<\infty\,.
\end{equation}

We now consider the second part of (\ref{aeqthm12}). Again by Potter's theorem and the fact that $\tilde{X}_{n-1}$ is independent of $W_{n-1}$ and $\Pi_{n-1}$ (under $\P$), we have
\begin{eqnarray}
\E\phi\left(\frac{|\tilde{X}_{n-1}|}{\Pi_{n-1}^{(\beta-1)/\beta}}\cdot W_{n-1}^{1/\beta}\right)
&=&\E\left(\frac{W_{n-1}}{\Pi_{n-1}^{\beta-1}}\right)^{\frac{\alpha}{\beta}}|\tilde{X}_{n-1}|^\alpha\ell\left(\left(\frac{W_{n-1}}{\Pi_{n-1}^{\beta-1}}\right)^{\frac{1}{\beta}}|\tilde{X}_{n-1}|\right)\nonumber\\
&\le&C\E\phi(|\tilde{X}_0|)\cdot( I_2^+(n)+I_2^-(n)),
\end{eqnarray}
where
\begin{eqnarray*}
I_2^+(n)=\E\left(\frac{W_{n-1}}{\Pi_{n-1}^{\beta-1}}\right)^{\frac{\alpha+\delta}{\beta}}, \\
I_2^-(n)=\E\left(\frac{W_{n-1}}{\Pi_{n-1}^{\beta-1}}\right)^{\frac{\alpha-\delta}{\beta}}.
\end{eqnarray*}
Here we have used the fact that under $\P$, each $\tilde{X}_{n-1}$ has the same distribution as $\tilde{X}_0$; $C=C(\ell_1,\delta,\beta)>0$ is a constant depending only on $\ell_1$, $\delta$ and $\beta$; $\delta\le \delta_0$.\\
We can estimate $I_2^+(n)$ as we have done for $I_1^+(n)$: we have
\begin{eqnarray}
I_2^+(n)&=&\E\left(\Pi_{n-1}^{-(\alpha+\delta)}\cdot\E_\zeta Z_{n-1}^{\frac{\alpha+\delta}{\beta}}\right)\nonumber\\
&\leq& \E\left(\Pi_{n-1}^{-(\alpha+\delta-\frac{\beta-1}{\beta})}\cdot\left(\E_\zeta Z_{n-1}^{\alpha+\delta-(\beta-1)}\right)^{\frac{1}{\beta}}\left(\E_{\zeta}W_{n-1}\right)^{\frac{\beta-1}{\beta}}\right)\nonumber\\
&=& \E\left(\Pi_{n-1}^{-(\alpha+\delta-\frac{\beta-1}{\beta})}\cdot\left(\E_\zeta Z_{n-1}^{\alpha+\delta-(\beta-1)}\right)^{\frac{1}{\beta}}\right)\nonumber\\
&\leq& \E\left(\Pi_{n-1}^{-(\alpha+\delta-\frac{\beta-1}{\beta})}\left(\E_\zeta Z_{n-1}^{\alpha+\delta-(\beta-1)+(\beta-1-2\delta)}\right)^{\frac{1}{\beta}}\right)\nonumber\\
&=&\E\left(\Pi_{n-1}^{-\frac{(\alpha+\delta-1)(\beta-1)}{\beta}}\cdot\left(\E_\zeta W_{n-1}^{\alpha-\delta}\right)^{\frac{1}{\beta}}\right)\nonumber\\
&\leq&\left(\E W_{n-1}^{\alpha-\delta}\right)^{\frac{1}{\beta}}\left(\E \Pi_{n-1}^{1-(\alpha+\delta)}\right)^{\frac{\beta-1}{\beta}}.
\end{eqnarray}
It follows that
\begin{equation}
\sum_{n\geq1}^\infty I_2^+(n)\leq\left(\sup_{n\geq1}\E W_{n-1}^{\alpha-\delta}\right)^{\frac{1}{\beta}}\cdot\left(\sum_{n=1}^\infty\left(\E m_0^{1-(\alpha+\delta)}\right)^{\frac{n(\beta-1)}{\beta}}\right)<\infty\,.
\end{equation}
Similarly we obtain
\begin{equation}
I_2^-(n)\le (\E W_{n-1}^{\alpha-\delta})^{\frac{1}{\beta}}\left(\E \Pi_{n-1}^{1-(\alpha-\delta)}\right)^{\frac{\beta-1}{\beta}}
\end{equation}
and
\begin{equation}
\sum_{n\geq1}^\infty I_2^-(n)\leq\left(\sup_{n\geq1}\E W_{n-1}^{\alpha-\delta}\right)^{\frac{1}{\beta}}\cdot\left(\sum_{n=1}^\infty\left(\E m_0^{1-(\alpha-\delta)}\right)^{\frac{n(\beta-1)}{\beta}}\right)<\infty\,.
\end{equation}
Therefore,
\begin{equation}\label{aeqthm111}
\sum_{n=1}^\infty\E\phi\left(\frac{|\tilde{X}_{n-1}|}{\Pi_{n-1}^{(\beta-1)/\beta}}W_{n-1}^{1/\beta}\right)<\infty.
\end{equation}
Combining (\ref{aeqthm12}), (\ref{aeqthm110}) and (\ref{aeqthm111}), we get
\begin{equation}\label{aeqthm112}
\E\phi(W^*-1)<\infty,
\end{equation}
which is equivalent to $\E\phi(W^*)<\infty$.

(ii) We next show that (b) implies (c). Obviously,
\begin{equation*}
\E\phi(W)\leq\E\phi(W^*)<\infty\,;
\end{equation*}
by Jensen's inequality, for any $n\ge1$,
\begin{equation*}
\E\phi(W_n)\ge\phi(\E W_n)=\phi(1)>0.
\end{equation*}
So by the dominated convergence theorem, we see that
\begin{equation*}
\E\phi(W)=\lim_{n\rightarrow\infty}\E\phi(W_n)\ge\phi(1)>0.
\end{equation*}

(iii) We finally show that (c) implies (a).
Notice that the limit $W$ satisfies the equation
\begin{equation}
W=\sum_{i=1}^{Z_1}\frac{W^{(i)}}{m_0},
\end{equation}
where under $\P_{\zeta}$, $(W^{(i)})$ are independent of each other, and have the same law as $W$ under $\P_{T\zeta}$: $\P_{\zeta}(W^{(i)}\in\cdot)=\P_{T\zeta}(W\in\cdot)$, $T$ being the usual translation: $T\zeta=(\zeta_1,\zeta_2,\ldots)$ if $\zeta=(\zeta_0,\zeta_1,\ldots)$. By Jensen's inequality, writting $\E_{\zeta,1}(\cdot)=\E_{\zeta}(\cdot|\F_1)$, we have
\begin{eqnarray}
\E_{\zeta}\phi(W)=\E_{\zeta}\phi\left(\sum_{i=1}^{Z_1}\frac{W^{(i)}}{m_0}\right)
&\ge& \E_{\zeta}\phi\left(\E_{\zeta,1} \sum_{i=1}^{Z_1}\frac{W^{(i)}}{m_0}\right)\nonumber\\
&=&\E_{\zeta}\phi(Z_1/m_0)=\E_{\zeta}\phi(W_1).
\end{eqnarray}
Therefore
\begin{eqnarray}
\E\phi(W_1)\le\E\phi(W).
\end{eqnarray}

\end{proof}


\noindent {\bf Remark.} \quad
For technical reasons, in the proof of Theorem \ref{athm1}, we have used the martingale $(W_n, \F_n)$ under $\P$ rather than the more natural martingale $(W_n, \mathcal{E}_n)$ under $\P_{\zeta}$. In fact, if we take the later martingale, then instead of (\ref{aeqlm13}), we have
\begin{equation}
\E_{\zeta}\phi(W^*-1)\le B\left(\E_{\zeta}\phi\left(\left(\sum_{n=1}^{\infty}\E_{\zeta}(|D_n|^{\beta}|\EE_{n-1})\right)^{\frac{1}{\beta}}\right)+\sum_{n=1}^{\infty}\E_{\zeta}\phi(|D_n|)\right);
\end{equation}
instead of (\ref{aeqlm14}), we obtain
\begin{equation}
\E_{\zeta}(|D_n|^{\beta}|\EE_{n-1})\le B\frac{\sigma_{n-1}^{\beta}(\beta)}{\Pi_{n-1}^{\beta-1}}W_{n-1}.
\end{equation}
Taking expectations and using the same argument as in part (i) of the proof of Theorem \ref{athm1}, we obtain 
\begin{equation}
\E\phi(W^*-1)\le C\left(\sum_{n=1}^{\infty}(\tilde{I}_1^+(n)+\tilde{I}_1^-(n))+\sum_{n=1}^{\infty}\E\phi(|D_n|)\right),
\end{equation}
where
\begin{eqnarray*}
\tilde{I}_1^+(n)=\E\Pi_{n-1}^{1-\beta}\A^{\frac{\alpha+\delta-\beta}{\beta}}[\sigma_{n-1}(\beta)]^{\alpha+\delta}W_{n-1}^{\frac{\alpha+\delta}{\beta}},\\
\tilde{I}_1^-(n)=\E\Pi_{n-1}^{1-\beta}\A^{\frac{\alpha-\delta-\beta}{\beta}}[\sigma_{n-1}(\beta)]^{\alpha-\delta}W_{n-1}^{\frac{\alpha-\delta}{\beta}}.
\end{eqnarray*}
The problem here is that we have to deal with the extra term $\sigma_{n-1}(\beta)$ in $\tilde{I}_1^{\pm}(n)$. We can do this by H\"older's inequality, but we then need an extra moment condition on $\sigma_{n-1}(\beta)$. Elementary calculations show that if for some positive number $\delta_0$, either (a) $\alpha<2$ and $\E[\sigma_0(\alpha)]^{\alpha(\alpha+\delta_0)}<\infty$, or (b) $\alpha\ge 2$ and $\E[\sigma_0(2)]^{2(\alpha+\delta_0)}<\infty$, then $\sum_{n=1}^{\infty}\tilde{I}_1^{\pm}(n)<\infty$, provided that $\E W_1^{\alpha}\ell(W_1)<\infty$. This leads a proof of Theorem \ref{athm1} under the preceding extra moment condition.



{}

\section{Proofs of Theorem \protect\ref{athm2} and Corollary \protect\ref%
{acor2}}

\setcounter{equation}{0}

\begin{proof}[Proof of  Theorem \ref{athm2}]
Write $\phi(x)=x\ell(x)$. By Lemma \ref{alm3}, we can assume that $\phi$ is convex on $(0,\infty)$, $\phi(x^{1/2})$ is concave on $(0,\infty)$ and $\ell$ is also concave on $(0,\infty)$ with $\ell(0)=0$.

Notice that under $\P$, $\tilde{X}_{n-1}$ is independent of $\Pi_{n-1}$. As $\ell$ is concave, we have
\begin{eqnarray}\label{aeqthm21}
\E \Pi_{n-1}\phi\left(\frac{|\tilde{X}_{n-1}|}{\Pi_{n-1}}\right)
&=& \E  |\tilde{X}_{n-1}|\ell\left(\frac{|\tilde{X}_{n-1}|}{\Pi_{n-1}}\right)\nonumber\\
&\le& \E |\tilde{X}_{n-1}|\ell(b^n|\tilde{X}_{n-1}|)\nonumber\\
&=&\E |\tilde{X}|\ell(b^n|\tilde{X}|),
\end{eqnarray}
where $\tilde{X}$ is a random variable having the same distribution as $(\tilde{X}_n)_{n\ge0}$ and $b=\E m_0^{-1}<1$. According to the inequality (\ref{aeqlm12}), we have
\begin{eqnarray}\label{aeqthm22}
\E\phi(W^*-1)
&\le& C\sum_{n=1}^{\infty}\E \Pi_{n-1}\phi\left(\frac{|\tilde{X}_{n-1}|}{\Pi_{n-1}}\right)\nonumber\\
&\le& C\sum_{n=1}^{\infty}\E |\tilde{X}|\ell(b^n|\tilde{X}|)\nonumber\\
&\le& C \sum_{n=1}^{\infty}\E |\tilde{X}|\int_{b^n|\tilde{X}|}^{b^{n-1}|\tilde{X}|}\frac{\ell(t)}{t}dt\nonumber\\
&=&C\E |\tilde{X}|\int_0^{|\tilde{X}|}\frac{\ell(t)}{t}dt\nonumber\\
&\le& C\E|\tilde{X}|(1+\hat{\ell}(|\tilde{X}|))<\infty\,.
\end{eqnarray}
This yields $\E W^*\ell(W^*)<\infty$, and
\begin{gather*}
\E W\ell(W)\le \E W^*\ell(W^*)<\infty\,.
\end{gather*}
If in addition, $\ell$ is slowly varying at $\infty$, then we can use Potter's theorem to replace the Jensen's inequality in (\ref{aeqthm21}), to relax the assumption $\E m_0^{-1}<1$.
Recall that for this $\ell$, we have shown that
\begin{eqnarray}
\E\phi(W^*-1)
&\le& C\sum_{n=1}^{\infty}\E |\tilde{X}|\ell\left(\frac{|\tilde{X}|}{\Pi_{n-1}}\right)\label{aeqthm23a}\\
&=& C\sum_{n=1}^{\infty}(I_3(n)+I_3'(n)),\label{aeqthm23}
\end{eqnarray}
where
\begin{eqnarray*}
I_3(n)=E|\tilde{X}|\ell\left(\frac{|\tilde{X}|}{\Pi_{n-1}}\right)\mathbf{1}_{\{\Pi_{n-1}^{-1}\le a^{n-1}\}},\\
I_3'(n)=E|\tilde{X}|\ell\left(\frac{|\tilde{X}|}{\Pi_{n-1}}\right)\mathbf{1}_{\{\Pi_{n-1}^{-1}> a^{n-1}\}},
\end{eqnarray*}
$a\in(0,1)$ will be determined later.
By the same argument as above, we get
\begin{equation}
\sum_{n=1}^{\infty}I_3(n)\le C\E|\tilde{X}|(\hat{\ell}(|\tilde{X}|)+1)<\infty\,.
\end{equation}
We now estimate $I_3'(n)$. For fixed $n$, we divide it into two parts:
\begin{eqnarray*}
I_{3,1}'(n)=\E|\tilde{X}|\ell\left(\frac{|\tilde{X}|}{\Pi_{n-1}}\right)\mathbf{1}_{\{\Pi_{n-1}^{-1}>a^{n-1}\}}\mathbf{1}_{\{|\tilde{X}|a^{n-1}>1\}};\\
I_{3,2}'(n)=\E|\tilde{X}|\ell\left(\frac{|\tilde{X}|}{\Pi_{n-1}}\right)\mathbf{1}_{\{\Pi_{n-1}^{-1}>a^{n-1}\}}\mathbf{1}_{\{|\tilde{X}|a^{n-1}\le 1\}}.
\end{eqnarray*}
As $\ell$ is increasing and slowly varying at $\infty$, by Potter's theorem, we have: for $\delta>0$,
\begin{eqnarray}\label{aeqthm24}
I_{3,1}'(n)
&\le& C\E|\tilde{X}|\ell(|\tilde{X}|a^{n-1})(\Pi_{n-1}a^{n-1})^{-\delta}\nonumber\\
&\le& C\E|\tilde{X}|\ell(|\tilde{X}|)(\Pi_{n-1}a^{n-1})^{-\delta}\nonumber\\
&=& C\E|\tilde{X}|\ell(|\tilde{X}|)\cdot(\E m_0^{-\delta}\cdot a^{-\delta})^na^{\delta}.
\end{eqnarray}
Let $\rho(x)=\E m_0^{-x}$. Since $\rho(\delta_0)<\infty$ and $\rho(x)$ is convex on $(0,\delta_0)$ with $\rho(0)=1$ and $\rho'(0)=-\E\ln m_0<0$, there exists some $\gamma_0>0$ such that
\begin{equation*}
\E m_0^{-x}<1, \qquad \forall x\in(0,\gamma_0).
\end{equation*}
Choose $\delta\in(0,\gamma_0)$, and let $0<a<1$ be defined by  $\E m_0^{-\delta}=a^{2\delta}$. Notice that $\E |\tilde{X}|\ell(|\tilde{X}|)\le C\E|\tilde{X}|(\hat{\ell}(|\tilde{X}|)+1)<\infty$. Therefore,
\begin{equation}
\sum_{n=1}^{\infty}I_{3,1}'(n)\le C\E|\tilde{X}|\ell(|\tilde{X}|)\cdot\sum_{n=1}^{\infty}a^{\delta(n+1)}<\infty\,.
\end{equation}
Similarly, using Potter's theorem in $I_{3,2}'(n)$, we get
\begin{eqnarray}
I_{3,2}'(n)
&\le& \E|\tilde{X}|\ell(\Pi_{n-1}^{-1}a^{1-n})\mathbf{1}_{\{\Pi_{n-1}^{-1}>a^{n-1}\}}\nonumber\\
&\le& C\E|\tilde{X}|\ell(1)(\Pi_{n-1}a^{n-1})^{-\delta}\nonumber\\
&\le& C\E|\tilde{X}| \cdot (\E m_0^{-\delta}a^{-\delta})^na^{\delta}\nonumber\\
&\le& C\E|\tilde{X}|\cdot a^{\delta(n+1)}.
\end{eqnarray}
Hence
\begin{equation}
\sum_{n=1}^{\infty}I_{3,2}'(n)\le C\E|\tilde{X}|\cdot \sum_{n=1}^{\infty}a^{\delta(n-1)}<\infty\,.
\end{equation}
Therefore, we have shown that
\begin{equation}
\E\phi(W^*-1)<\infty\,,
\end{equation}
which is equivalent to $\E\phi(W^*)<\infty$.

\end{proof}

\begin{proof}[Proof of Corollary \ref{acor2}]
Let
\begin{equation}
\phi_1(x)=
\begin{cases}
\frac{\phi'(1)}{2}x^2 &\qquad {\rm if}\ x\le 1;\\
\phi(x)+c_0 &\qquad {\rm if}\ x>1
\end{cases}
\end{equation}
where $\phi(1)+c_0=\frac{\phi'(1)}{2}$. Then it is easily seen that $\phi_1\asymp\phi, \phi_1(0)=0,\phi_1'(0+)=0$ and $\int_0^1\frac{\phi_1'(t)}{t}dt=\phi'(1)<\infty$. Moreover, $\phi_1$ is convex with positive concave derivative $\phi_1'$ on $(0,\infty)$, so that the function $x\mapsto \phi_1(x^{1/2})$ is concave on $(0,\infty)$. Applying the BDG-inequality and the concavity of $\phi_1(x^{1/2})$ (which implies the subadditivity), we obtain
\begin{eqnarray}\label{aeqcor21}
\E\phi_1(W^*-1)
&\le& C\E\phi_1\left(\left(\sum_{n=1}^{\infty}|D_n|^2\right)^{\frac{1}{2}}\right)\nonumber\\
&\le& C\sum_{n=1}^{\infty}\E\phi_1(|D_n|),
\end{eqnarray}
where $C=C(\phi_1)>0$ is a constant depending only on $\phi_1$.

Recalling that under $\P_{n-1}$, $D_n$ is a sum of a sequence of martingale differences with respect to $(\tilde{\F}_k)$. Hence, again by the BDG-inequality applied to $D_n$, and the concavity of $\phi_1(x^{1/2})$, we get
\begin{eqnarray}\label{aeqcor22}
\E_{n-1}\phi_1(|D_n|)
&\le& C\E_{n-1}\phi_1\left(\left(\sum_{|u|=n-1}\frac{|\tilde{X}|^2}{\Pi_{n-1}^2}\right)^{\frac{1}{2}}\right)\nonumber\\
&\le& C\E_{n-1}\sum_{|u|=n-1}\phi_1\left(\frac{|\tilde{X}|}{\Pi_{n-1}}\right)\nonumber\\
&=& C Z_{n-1}\cdot\E_{n-1}\phi_1\left(\frac{|\tilde{X}_{n-1}|}{\Pi_{n-1}}\right)
\end{eqnarray}
where $C>0$ is independent of $n$.
Taking integral on both sides of the inequality above, and noting that $\phi_1'$ is concave, we obtain:
\begin{eqnarray}\label{aeqcor23}
\E\phi_1(|D_n|)
&\le& C\E\E_{\zeta}\left(Z_{n-1}\cdot\E_{n-1}\phi_1\left(\frac{|\tilde{X}_{n-1}|}{\Pi_{n-1}}\right)\right)\nonumber\\
&=&C\E\Pi_{n-1}\phi_1\left(\frac{|\tilde{X}_{n-1}|}{\Pi_{n-1}}\right)\nonumber\\
&=&C\E|\tilde{X}_{n-1}|\int_0^1\phi_1'\left(\frac{|\tilde{X}_{n-1}|}{\Pi_{n-1}}s\right)ds\nonumber\\
&\le& C\E|\tilde{X}_{n-1}|\phi_1'(b^n|\tilde{X}_{n-1}|)\nonumber\\
&=&C\E |\tilde{X}|\phi_1'(b^n|\tilde{X}|),
\end{eqnarray}
where $\tilde{X}$ is a random variable having the same distribution as $(\tilde{X}_n)_{n\ge0}$ and $b=\E m_0^{-1}<1$.
Similarly to (\ref{aeqthm22}), combining (\ref{aeqcor21}) and (\ref{aeqcor23}), we obtain
\begin{eqnarray}\label{aeqcor24}
\E\phi_1(W^*-1)
&\le& C\E\sum_{n=1}^{\infty}|\tilde{X}|\phi_1'(b^n|\tilde{X}|)\nonumber\\
&\le& C\E|\tilde{X}|\int_0^{|\tilde{X}|}\frac{\phi_1'(t)}{t}dt\nonumber\\
&\le& C\E|\tilde{X}|(\tilde{\phi}_1(|\tilde{X}|)+1).
\end{eqnarray}
As $\phi\asymp\phi_1$ and $\tilde{\phi}\asymp\tilde{\phi}_1$, this yields
\begin{equation}\label{aeqcor25}
\E\phi(W^*-1)\le C\E |\tilde{X}|(\tilde{\phi}(|\tilde{X}|)+1)<\infty\,.
\end{equation}
Therefore $\E\phi(W^*)<\infty$,
and
\begin{equation*}
\E\phi(W)\le \E\phi(W^*)<\infty\,.
\end{equation*}

\end{proof}


\section{Proof of Theorem \protect\ref{athm3}}
\setcounter{equation}{0}
Before giving the proof of Theorem \ref{athm3}, we first show an extension of a theorem of Hsu and Robbins (1947) (see also Erd\"os (1949) or Baum and Katz (1965)). As usual, for a random variable $X$, we write $X^+=\max(X,0)$ and $X^-=\max(-X,0)$.
\begin{lemma}\label{alm4}
Let $(X_i)$ be i.i.d. with $m=\E X_1\in[-\infty,\infty)$. If $\E(X_1^+)^2<\infty$, then for all $a>m$
\begin{equation}\label{aeqlm41}
\sum_{n=1}^{\infty}\P\left(\sum_{i=1}^n X_i>na\right)<\infty\,.
\end{equation}
\end{lemma}

\begin{proof}
The result is due to Hsu and Robbins if $\E X_1^2<\infty$. What is new here that we may have $\E X_1^-=\infty$ or $\E (X_1^-)^2=\infty$. Notice that $\forall a\in\mathbb{R}$,
\begin{equation}\label{aeqlm42}
\P\left(\sum_{i=1}^n X_i>na\right)\le\P\left(\sum_{i=1}^n X_i^+>na_1\right)+\P\left(-\sum_{i=1}^n X_i^->na_2\right)
\end{equation}
where $a_1+a_2=a$. By the theorem of Hsu and Robbins (1947),
\begin{equation}\label{aeqlm43}
\sum_{n=1}^{\infty}\P\left(\sum_{i=1}^n X_i^+>na_1\right)<\infty\qquad \forall\ a_1>\E X_1^+.
\end{equation}
Now for all $C>0$,
\begin{equation*}
\P\left(-\sum_{i=1}^n X_i^->na_2\right)\le\P\left(\sum_{i=1}^n\max(-X_i^-, -C)>na_2\right).
\end{equation*}
Therefore, again by the theorem of Hsu and Robbins,
\begin{eqnarray}\label{aeqlm44}
\sum_{n=1}^{\infty}\P\left(-\sum_{i=1}^n X_i^->na_2\right)<\infty\qquad \forall\ a_2
&>&\E\max(-X_1^-,-C)\nonumber\\
&=&-\E\min(X_1^-,C).
\end{eqnarray}
Notice that $\lim\limits_{C\rightarrow\infty}\E\min(X_1^-,C)=\E X_1^-$ by the monotone convergence theorem. Hence (\ref{aeqlm44}) holds for all $a_2>-\E X_1^-(\ge -\infty)$. It follows from (\ref{aeqlm42}), (\ref{aeqlm43}) and (\ref{aeqlm44}) that (\ref{aeqlm41}) holds for all $a>\E X_1^+-\E X_1^-=\E X_1$.
\end{proof}


\vskip 2mm
\begin{proof}[Proof of Theorem \ref{athm3}]
Let
\begin{equation*}
\ell(x)=
\begin{cases}
1-\frac{1}{2x}, &\qquad {\rm if}\quad x>1;\\
\frac{x}{2},&\qquad{\rm if}\quad x\le1.
\end{cases}
\end{equation*}
Then $\phi(x)=x\ell(x)$ is convex, and the function $x\mapsto\phi(x^{1/2})$ is concave. By an argument similar to that in the proof of Theorem \ref{athm2}, we get (cf. (\ref{aeqthm23a}))
\begin{equation}
\E\phi(W^*-1)\le C \sum_{n=1}^{\infty}\E|\tilde{X}|\ell\left(\frac{|\tilde{X}|}{\Pi_{n-1}}\right).
\end{equation}

Let $b\in(e^{-\E\ln m_0},1)$ (by convention $e^{-\E\ln m_0}=0$ if $\E \ln m_0=+\infty$). For $n\ge0$, we divide the domain of integration above into two parts according to $\{\Pi_n^{-1}\le b^n\}$ or $\{\Pi_n^{-1}>b^n\}$, so that
\begin{equation}
\E\phi(W^*-1)\le C\sum_{n=0}^{\infty}(I_4(n)+I_4'(n)),
\end{equation}
where
\begin{eqnarray*}
I_4(n)=\E|\tilde{X}|\ell(|\tilde{X}|\Pi_n^{-1})\mathbf{1}_{\{\Pi_n^{-1}\le b^n\}},\\
I_4'(n)=\E|\tilde{X}|\ell(|\tilde{X}|\Pi_n^{-1})\mathbf{1}_{\{\Pi_n^{-1}> b^n\}}.
\end{eqnarray*}

We first estimate $I_4(n)$. Noting that $\ell$ is increasing on $[0,\infty)$, we get $I_4(n)\le\E|\tilde{X}|\ell(|\tilde{X}|b^n)$; moreover,
\begin{eqnarray}
\sum_{n=0}^{\infty}I_4(n)
&\le& C\E|\tilde{X}|\int_0^{|\tilde{X}|}\frac{\ell(t)}{t}dt\nonumber\\
&\le& C\E|\tilde{X}|(1+\ln^+|\tilde{X}|)<\infty\,.
\end{eqnarray}

To estimate $I_4'(n)$, as $\ell$ is bounded by $1$, we have
\begin{eqnarray}\label{*}
\sum_{n=0}^{\infty}I_4'(n)
&\le& \E|\tilde{X}|\cdot\sum_{n=0}^{\infty}\E\mathbf{1}_{\{\Pi_n^{-1}>b^n\}}\nonumber\\
&=& \E|\tilde{X}|\cdot\sum_{n=0}^{\infty}\P(\Pi_n^{-1}>b^n).
\end{eqnarray}
By Lemma \ref{alm4}, the sum on the right side of (\ref{*}) is finite if $\E\left(\ln^+\frac{1}{m_0}\right)^2<\infty$. Therefore,
\begin{equation}
\E\phi(W^*)<\infty\,,
\end{equation}
which is equivalent to $\E W^*<\infty$.

\end{proof}


\end{document}